\documentclass[12pt]{elsarticle}
\usepackage{amsfonts}
\usepackage{amsmath}
\usepackage{graphicx}
\usepackage{epsfig}
\usepackage{amssymb}
\usepackage{amstext}
\usepackage{tkz-berge}
\newtheorem{theorem}{Theorem}[section]
\newtheorem{lemma}{Lemma}[section]

\newenvironment{PrfFact}{{\bf Proof }}{{\hfill\tiny{$\blacksquare$\\}}}
\begin{document}
\begin{frontmatter}

\title{Remarks on Barnette's Conjecture}

\author{Jan~Florek}
\ead{jan.florek@pwr.edu.pl}

\address{Faculty of Pure and Applied Mathematics,

 Wroclaw University of Science and Technology,

 Wybrze\.{z}e Wyspia\'nskiego 27,
50-370 Wroc{\l}aw, Poland}

\fntext[]{This research did not receive any specific grant from funding agencies in the public, commercial, or not-for-profit sectors.}

\begin{abstract}%Barnette's conjecture is the statement that every cubic $3$-connected  bipartite planar graph is hamiltonian. 
Let $P$ be a cubic $3$-connected  bipartite plane graph which has a $2$-factor which consists only of facial $4$-cycles, and suppose that $P^{*}$ is the dual graph.  We show that $P$ has at least  $3^{\frac{2|P^{*}|}{\Delta^{2}{(P^{*})}}}$ different Hamilton  cycles.

%\footnotetext{2010 \textit{Mathematics Subject Classification}: 05C45, 05C10.}
%\footnotetext{\textit{Key words and phrases}: Barnette's conjecture, Hamilton cycle, induced tree, contraction of a path.}

 \end{abstract}


\end{frontmatter}
\section{Introduction}

%All graphs considered in this paper are finite and simple. We use \cite{flobar1} as reference for undefined terms. In particular, $ V(G)$ denotes  

%In 1884, P.~G.~Tait posed a conjecture that all cubic $3$-connected plane graphs had a hamiltonian cycle. The conjecture was of a special interest since it implied the four color theorem. Because of counterexamples found the conjecture has been modified in various manners. In particular,

Let ${\cal P}$ denote the family of cubic $3$-connected bipartite plane graphs. Barnette, in 1969  (\cite{flobar8}, Problem~5), conjectured that every graph in ${\cal P}$ has a Hamilton cycle. In \cite{flobar4}, Goodey proved that if a graph in ${\cal P}$ has only faces with $4$ or $6$ sides, then it is hamiltonian (see also Feder and Subi \cite{flobar3}, and Bagheri, Feder, Fleischner and Subi \cite{flobar0}).
 Holton, Manvel and McKay \cite{flobar5} used computer search to confirm Barnette's conjecture for graphs up to $64$ vertices. The problem whether a cubic bipartite planar graph has a Hamilton cycle (without the assumption of $3$-connectivity) is NP-complete, as shown by Takanori, Takao, and Nobuji \cite{flobar7}.

All graphs considered in this paper are finite and simple. We use \cite{flobar1} as reference for undefined terms. In particular, if $P$ is a plane graph, then $ V(P)$  is its vertex set and $E(P)$ is its edge set.  The set of neighbours of a vertex $v$ in $P$ is denoted by $N(v)$. $|P|$ is the number of vertices of $P$ and $\Delta (P)$ is the maximum degree of $P$. Hence, if $P^{*}$ is the dual graph of $P$, then $|P^{*}|$  is the number of faces in $P$, and $\Delta (P^{*})$  is the maximum order of all face boundaries in $P$.  

Let  ${\cal P}(4)$  be the family of all graphs in ${\cal P}$ which possess a $2$-factor  consisting  only of facial $4$-cycles. The  following theorem yields a  lover bound on the number of hamiltonian cycles in every graph belonging to ${\cal P}(4)$.

\begin{theorem}\label{theorem1.1} Every graph $P \in {\cal P}(4)$ has at least $3^{\frac{2|P^{*}|}{\Delta^{2}{(P^{*})} } }$ different Hamilton cycles.
\end{theorem}

A \textit{sequence of  faces} in a plane graph $P$ is a sequence  of different faces  in $P$ such that two faces from the sequence are adjacent if and only if they are consecutive in the sequence. A sequence of faces $f_{0}f_{1} \ldots f_{k}$  in $P$ is called a sequence from $f_0$ to $f_k$ with \textit{length}  $k$.  The distance  of two faces $f, h$ in $P$ is the length of a shortest sequence of faces from $f$ to $h$. We prove the following theorem which generalizes results obtained by Florek in \cite{flobar2}.
\begin{theorem}\label{theorem1.2} Let $P\in {\cal P}(4)$ and suppose that $M$ is a set of faces  in P  such that the distance of any two  different faces in $M$ is greater than $4$.  If for every face in $M$ any edge is chosen on this face, there is a Hamiltonian cycle containing all other edges of faces in $M$.
\end{theorem}

We will use a result that expresses hamiltonicity in the terms of the dual graph.  Let ${\cal E}$ be the dual family to ${\cal P}$. Hence, ${\cal E}$ is the family of all Eulerian plane triangulations. Let $G \in {\cal E}$.  Stein \cite{flobar6} proved that $G^{*}$ is hamiltonian if and only if $G$ possesses two disjoint induced tree-subgraphs which together contain all  vertices of $G$ (if $X$, $Y$ are such tree-subgraphs in $G$, then the set $ \{e^* \in E(G^*) : e \hbox { is an } X-Y \hbox { edge} \}$ is  the edge set of a Hamilton cycle in $G^*$).

Let ${\cal E}(4)$ be the dual family to ${\cal P}(4)$.  Our aim is to prove the following two theorems which express the above two theorems in the dual form:
\begin{theorem}\label{theorem1.3}
Every graph $G \in {\cal E}(4)$ has at least  $3^{\frac{2|G|}{\Delta^{2}{(G)}}}$ different pairs of disjoint  induced acyclic subgraphs which together contain all the vertices of~$G$.
\end{theorem}
\begin{theorem}\label{theorem1.4} 
Let $G \in {\cal E}(4)$ and suppose that $L \subset V(G)$ is a set of vertices such that the distance of every two different vertices of $L$ is greater than $4$.  If for every $v \in L$ any vertex $n(v) \in N(v)$ is chosen, then there are two induced tree-subgraphs which together contain all vertices of $G$ and one of them contains the set $N(v)   \backslash \{n(v)\}$, for every $v \in L$.
\end{theorem}
 
\section{Main results}
Let $G\in {\cal E}(4)$. It is well known that every Eulerian plane triangulation has a proper $3$-colouring of vertices. Let $\{B, W, R\}$ be a vertex partition of $V(G)$ into independent  black, white  and red colour classes. Certainly, one of them (say $R$) contains only vertices of degree $4$.

Florek \cite{flobar2} proved that every graph in ${\cal E}(4)$ is $4$-connected. Hence, the set of neighours of any vertex in $G$ induces a cycle of $G$. For a path $svt$, which is disjoint with $R$, we define an \textsl{operation} $\alpha = \alpha(s,v,t)$ in the following way: let $v_{1} \ldots  v_{n} v_{1}$ be a cycle induced by neighbours of $v$ ($s = v_1$ and $t = v_k$). Replace the vertex $v$ with a path $xuy$ and join the vertices of this path with the former neighbours of $v$ -- provided that $x$ is adjacent to $v_{1}, \ldots, v_{k}$,  $u$ is adjacent to $v_1$ and $v _k$, and $y$  is adjacent to $v_{k}, \ldots,  v_{n}, v_{1}$. Let $\alpha(G)$ be the graph constructed from $G$ by $\alpha$. We put $\alpha(R) = R\cup \{u\}$, $\alpha(W) = (W\backslash \{v\})\cup \{x, y\}$ and  $\alpha(B) = B$, for $v \in W$ (similarly, $\alpha(B) = (B\backslash \{v\})\cup \{x, y\}$ and  $\alpha(W) = W$, for $v \in B$, respectively). Notice that  $\alpha(G) \in  {\cal E}(4)$, because $\{\alpha(B), \alpha(W), \alpha(R)\}$ is a vertex partition of $V(\alpha(G))$ into independent colour classes such that  $\alpha(R)$ contains only vertices of degree $4$.  It was proved in \cite{flobar2} that every graph $G\in {\cal E}(4)$ can be constructed (up to isomorphism) from the octahedron by iterating the operation $\alpha$. Hence, by induction we obtain:
\\
(i) $|B \cup W| >  \frac{1}{2}|G|$.

%We say that an induced subgraph $A$ of $G$ is \textit{black-closed} (\textit{white-closed}) if every black (white, respectively) vertex, which  is adjacent with a vertex belonging to $A$, belongs to $A$.

 We say that a pair $(C, D)$ of disjoint induced subgraphs of $G$ is \textit{(black,white)-closed} if the following two conditions are satisfied:

(1) if $v \in B$ ($v \in B\backslash V(D)$)  is adjacent with a white (red, respectively) vertex of $C$, then $v$ is vertex of $C$,

(2) if $v \in W$ ($v \in W\backslash V(C)$)  is adjacent with a black (red, respectively) vertex of $D$, then $v$ is vertex of $D$.

\begin{lemma}\label{lem2.1}
Let $G \in {\cal E}(4)$ and suppose that $V(G)$ has a vertex partition into independent red, black and white colour classes such that the red colour class contains only vertices of degree $4$.  Every  pair of disjoint (black, white)-closed  induced acyclic subgraphs of $G$ can be extended to such a pair of  disjoint induced acyclic subgraphs  which together contain all the vertices of $G$.
\end{lemma}
\begin{PrfFact}
Let $(C, D)$ be a pair of  disjoint (black,white)-closed  induced acyclic subgraphs of $G$.  First we add to $C$ (or $D$) all black  (white, respectively) vertices which not belong to $C\cup D$. Then we obtain a pair  $(C', D')$ of disjoint induced acyclic subgraphs of $G$. Next we consider any red vertex~$v$ not belonging to $C \cup D$.  Since $(C', D')$  is (black,white)-closed and all neighbours of $v$ belong to $C' \cup D'$  one of the following  cases occurs:
\\
(j) two black neighbours of $v$ belong to $C'$,
\\
(jj) two white neighbours of $v$ belong to $D'$.
\\
If  two black neighbours of $v$ belong to a path contained in $C'$, then we add $v$ to $D'$. If two white neighbours of $v$ belong to a path contained in $D'$ then we add $v$ to $C'$.  Otherwise, we add $v$ to $C'$ or $D'$.  Then we obtain a pair  $(X, Y)$ of disjoint  induced acyclic subgraphs of $G$ such that
$C' \subset X$, $D' \subset Y$ and $V(X) \cup V(Y) = V(G)$. 
\end{PrfFact}

Now we are ready to prove Theorem \ref{theorem1.3} and Theorem \ref{theorem1.4}
%\begin{theorem}\label{theo2.1}
%Every graph $G \in {\cal E}(4)$ has at least  $3^{\frac{2|G|}{\Delta^{2}{(G)}}}$ different pairs of disjoint  induced acyclic subgraphs which together contain all the vertices of~$G$.
%\end{theorem}
$$\mbox{\bf Proof  of Theorem \ref{theorem1.3}}$$
\begin{PrfFact}
%Let $\{B, W, R\}$ be a vertex partition of $V(G)$ into independent  black, white  and red colour classes such that the $R$ contains only vertices of degree~$4$. 
Suppose that $G[B\cup W]$ is a subgraph of $G$ induced by all black and white vertices of $G$. For every different vertices $a, b$ of $G[B\cup W]$ which have  the same colour and  $d_{G}(a,b) = 2$  we add the edge $ab$ to  $G[B\cup W]$. Then we obtain a graph  $J$ such that 
$$\Delta(J) \leqslant \frac{1}{2}\Delta(G) + \frac{1}{2}\Delta(G) (\frac{1}{2}\Delta(G)-1) - 1 < \frac{1}{4}\Delta^{2}(G).$$
Thus, by the greedy algorithm, $\chi(J) \leqslant \frac{1}{4}\Delta^{2}(G)$. Hence, $J$ has an independent set of vertices $K$ which, by (i), has at least
$$|K| \geqslant \frac{|J|}{\chi(J)} =  \frac{|B\cup W|}{\chi(J)} > \frac{2|G|}{\Delta^{2}(G)}$$
vertices.  

Suppose that for every $v \in K$  any red vertex $n(v) \in N(v)$ is chosen.   Since $G$ is $4$-conected  $N(v) \backslash \{r(v)\}$  induces a path-subgraph in $G$, denoted as $P(v,n(v))$.  If $v$ is black (white), then $(N(v) \cap W) \cup \{v\}$ ($(N(v) \cap B)\cup \{v\}$, respectively) induces a star-subgraph in $G$, denoted as $S(v)$.  

Let $\{B_1,B_2,W_1,W_2\}$ be a partition of the set $K$ such that vertices in $B_1\cup B_2$ are black and in $W_1 \cup W_2$ are white. Notice that 
$$C = \bigcup_{v \in W_1} P(v, n(v)) \cup \bigcup_{v \in W_2} S(v) \hbox{ and }
D =  \bigcup_{v \in B_1} P(v, n(v)) \cup \bigcup_{v \in B_2} S(v)$$
are disjoint  induced acyclic subgraphs of $G$ because $d_{G}(a,b)\geqslant 4$ ($\geqslant 3$),  for every different vertices  $a, b \in K$ of the same colour class (of different colour classes, respectively). Remark that every black (white) vertex, which  is adjacent with a vertex belonging to $C$ (or $D$) belongs to $C$ ($D$, respectively). Hence follows that $G$ has at least
$$\prod _{v \in K} (\frac{deg_{G}(v)}{2} + 1)  \geqslant  3^{|K|} \geqslant 3^{\frac{2|G|}{\Delta^{2}(G)}}$$
 different pairs of disjoint  (black,white)-closed induced acyclic subgraphs.  Thus, by Lemma \ref{lem2.1}, $G$ has at least  $3^{\frac{2|G|}{\Delta^{2}{(G)}}}$ different pairs of disjoint  induced acyclic subgraphs which together contain all the vertices of $G$.
\end{PrfFact}
%\begin{theorem}\label{theo2.2} 
%Let $G \in {\cal E}(4)$ and suppose that $L \subset V(G)$ is a set of vertices such that $d(a,b) \geqslant 5$ for   $a \neq b \in L$.  If for every $v \in L$ any vertex $n(v) \in N(v)$ is chosen, then there are two induced tree-subgraphs which together contain all vertices of $G$ and  one of them contains the set $N(v)   \backslash \{n(v)\}$, for every $v \in L$. 
%\end{theorem}
$$\mbox{\bf Proof  of Theorem \ref{theorem1.4}}$$
\begin{PrfFact}
%Let $\{B, W, R\}$ be a vertex partition of $V(G)$ into independent  black, white  and red colour classes such that the $R$ contains only vertices of degree~$4$. 
Let $L \subset V(G)$ be a set of vertices such that $d(a,b) \geqslant 5$ for every different vertices $a, b \in L$.

Suppose that for every $v \in V(G)$  any vertex $n(v) \in N(v)$ is chosen.  Since $G$ is $4$-conected  $N(v) \backslash \{n(v)\}$ induces a path-subgraph in $G$, denoted as $P(v, n(v))$. Let $R_b$ (or $R_w$) be the set of all red vertices such that $n(v)$ is black (white, respectively). If $v \in W \cup R_b$  ($v \in B \cup R_w$), then $(N(n(v)) \cap W) \cup \{n(v)\}$ ($(N(n(v)) \cap B)\cup \{n(v)\}$, respectively) induces a star-subgraph in $G$, denoted as $S(n(v))$. Notice that if $v$ is white (black), then the pair $(P(v, n(v)), S(n(v))$ ($(S(n(v)), P(v, n(v)))$, respectively) is (black,white)-closed.
Hence, 
$$(\bigcup_{v \in W \cap L} P(v,n(v)) \cup \bigcup_{v \in (B \cup R_w) \cap L} S(n(v)), \bigcup_{v \in B \cap L} P(v,n(v)) \cup \bigcup_{v \in (W \cup R_b)\cap L} S(n(v))$$
is a pair of disjoint (black,white)-closed  induced acyclic subgraphs of $G$, because $d(a,b)\geqslant 5$, for every different  $a, b \in L$. 

Hence, by Lemma \ref{lem2.1}, $G$ has a pair $(X, Y)$ of disjoint induced acyclic subgraphs  which together contain all the vertices of $G$ and such that $P(v,n(v)) \subset X$, for $v \in W\cap L$,  $S(n(v))\subset X$, for $v \in R_w\cap L$,  $P(v,n(v)) \subset Y$, for $v \in B\cap L$ and $S(n(v))\subset Y$, for $v \in R_b\cap L$.
\end{PrfFact}

\end{document}